\documentclass[12pt]{amsart}

\usepackage{tikz}
\usepackage{graphicx}
\usepackage{amsmath, amsfonts, amssymb}
\usepackage{mathtools}

\theoremstyle{theorem}
\newtheorem{theorem}{Theorem}
\newtheorem{proposition}[theorem]{Proposition}
\newtheorem{lemma}[theorem]{Lemma}
\newtheorem{corollary}[theorem]{Corollary}

\theoremstyle{definition}

\DeclareMathOperator{\C}{\mathcal{C}}
\DeclareMathOperator{\R}{\mathcal{R}}
\DeclareMathOperator{\Lf}{\mathcal{L}}
\DeclareMathOperator{\I}{\mathcal{I}}
\DeclareMathOperator{\Ou}{\mathcal{O}}

\begin{document}

\title{A More Malicious Maitre d'}
\markright{The spiteful table}

\author[R. Acton]{Reed Acton}
\author[T. K. Petersen]{T. Kyle Petersen}
\author[B. Shirman]{Blake Shirman}
\author[D. Toal]{Daniel Toal}
\address{Department of Mathematical Sciences, DePaul University, Chicago, IL}

\maketitle

\begin{abstract}
In this paper we study the problem of the Malicious Maitre d', as described in Peter Winkler's book \emph{Mathematical Puzzles: A Connoisseur's Collection}. This problem involves seating diners around a circular table with napkins placed between each pair of adjacent settings. The goal of the maitre d' is to seat the diners in a way that maximizes the number of diners who arrive at the table to find the napkins on both the left and right of their place already taken by their neighbors. Winkler proposes a solution to the problem that he claims is optimal. We analyze Winkler's solution using tools from enumerative combinatorics, then present a new strategy that performs better.
\end{abstract}

\section{The adaptive maitre d'}

We begin quoting directly from the book \textit{Mathematical Puzzles: A Connoisseur's Collection}, by Peter Winkler \cite[p.~22]{Winkler}:

\begin{quote}
\textsc{The Malicious Maitre D'}

{\small At a mathematics conference banquet, 48 male mathematicians,
none of them knowledgeable about table etiquette, find themselves
assigned to a big circular table. On the table, between each pair of
settings, is a coffee cup containing a cloth napkin. As each person is
seated (by the maitre d'), he takes a napkin from his left or right;
if both napkins are present, he chooses randomly (but the maitre d'
doesn't get to see which one he chose).

In what order should the seats be filled to maximize the expected
number of mathematicians who don't get napkins?

\ldots This problem can be traced to a particular event. Princeton
mathematician John H. Conway came to Bell Labs on March 30, 2001 to
give a ``General Research Colloquium." At lunchtime, [Winkler] found
himself sitting between Conway and computer scientist Rob Pike (now of
Google), and the napkins and coffee cups were as described in the
puzzle. Conway asked how many diners would be without napkins if they
were seated in \emph{random} order, and Pike said: ``Here's an easier
question---what's the \emph{worst} order?"  }
\end{quote}

There are two distinct problems here: Pike's problem, which Winkler calls ``The Malicious Maitre d','' and Conway's problem, which (later in the book) Winkler calls ``Napkins in a Random Setting.''

A solution to Conway's problem is given in Winkler's book, and the result is re-proved and generalized in several follow-up papers \cite{CP, Eriksen, Sudbury}. In the case that every diner reaches right or left with equal probability, we find the proportion of napkinless diners is about $12\%$. In fact the proportion of napkinless diners approaches $(2- \sqrt{e})^2 \approx .12339675$ as the number of diners tends to infinity.

This paper has its beginning in a reading course that examined the paper by Claesson and Petersen \cite{CP}, with the aim of reproducing the main result for Conway's napkin problem. In the introduction to that paper, the following is written about Pike's problem:

\begin{quote}
{\small The problem of the malicious maitre d' is not horribly difficult; if
you're having trouble finding a solution, you can see Winkler's book
for a nice explanation. }
\end{quote}

The authors of this paper read that line and were intrigued. How difficult is ``horribly'' difficult? We spent some time thinking about the problem, but since our main goal was to learn about Conway's problem, we didn't think for too long before reading the solution provided by Winkler. There, Winkler observes that since the maitre d' is not allowed to see which napkin a diner chooses, the strategy of the maitre d' is \emph{non-adaptive}. That is, the maitre d' is simply deciding ahead of time on a seating arrangement (i.e., a permutation) of the diners that will maximize the expected proportion of napkinless diners. Winkler then describes a permutation that achieves an expected proportion of $9/64$ (approximately $14\%$) of the diners to be napkinless. Winkler's justification that no permutation could do better was convincing to us. 

But this was not quite the end of the story. 
In Winkler's words \cite[p.~25]{Winkler}:
\begin{quote}
{\small If the maitre d' sees which napkin is grabbed each time he seats a diner (computer theorists would call him an ``adaptive adversary''), it is not hard to see that his best strategy is as follows. If the first diner takes (say) his right napkin, the next is seated two spaces to his right so that the diner in between may be trapped. If the second diner also takes his right napkin, the maitre d' tries again by skipping another chair to the right. If the second diner takes his left napkin (leaving the space between him and the first diner napkinless), the third diner is seated directly to the second diner's right. Further diners are seated according to the same rule until the circle is closed, then the remaining diners (some of whom are doomed to be napkinless) are seated. This results in 1/6 of the diners without napkins, on average.}
\end{quote}
In this paragraph Winkler both provides us with a new puzzle based on a different interpretation of Pike's question, and then presents a solution to that new puzzle. We will refer to this new puzzle as \textbf{The Adaptive Maitre d'}, and the rest of this paper is devoted to this puzzle.

\section{Trap setting versus napkin shunning}

We call Winkler's approach to the adaptive maitre d' puzzle \textbf{trap setting}. After some reflection, we were happy to agree with Winkler's conclusion that the proportion of napkinless diners is about $1/6$, as with $n$ seats there can be at most $n/3$ traps and about half of these will succeed. However, we had real trouble understanding the claim of optimality. Eventually we realized the reason we had a hard time proving trap setting is an optimal strategy was because it is \emph{not} an optimal strategy. 

In this paper we will present a better strategy: a \emph{more malicious maitre d'}!

The insight for the new strategy begins with the observation that there is a bijection between napkins at the table and diners at the table. Hence, for each napkinless diner, there must be, somewhere on the table, a dinerless napkin ignored by both of its adjacent diners. Rather than directly attempting to set traps by having diners sit apart and hoping they reach toward each other, the \textbf{napkin shunning} strategy seats a new diner on the side \emph{away} from the napkin the previous diner has chosen, in the hopes that the second diner reaches away from the first diner, shunning the napkin between. Once we succeed in shunning a napkin, the next diner is seated in the middle of a gap of empty seats and we attempt to shun another napkin. We continue to split gaps and try to shun napkins until all gaps filled. Since every shunned napkin will ultimately correspond to a napkinless diner, the traps will have set themselves.

Intuitively, the upside to napkin shunning is that when we fail to shun a napkin, only two seats at the table have been used, whereas a failure to set a trap in Winkler's algorithm uses three seats. As we will demonstrate in Section~\ref{sec:shun}, this extra available table space has a significant effect on the expected number of napkinless diners.

Our main result is the following.

\begin{theorem}[Napkin shunning is superior to trap setting]\label{thm:main}
The expected proportion of napkinless diners on a table with $n \geq 3$ seats:
\begin{enumerate}
\item is given by $\frac{(3n-2) -2(-1/2)^n}{18n} < 1/6$ when using the trap setting strategy, and
\item is at least $1/6=8/48$ and at most $3/16=9/48$ when using the napkin shunning strategy. Moreover, for all $n\geq 5$, this proportion is between $0.1769$ and $0.1831$.
\end{enumerate}
\end{theorem}

\begin{figure}
\includegraphics[width=12cm]{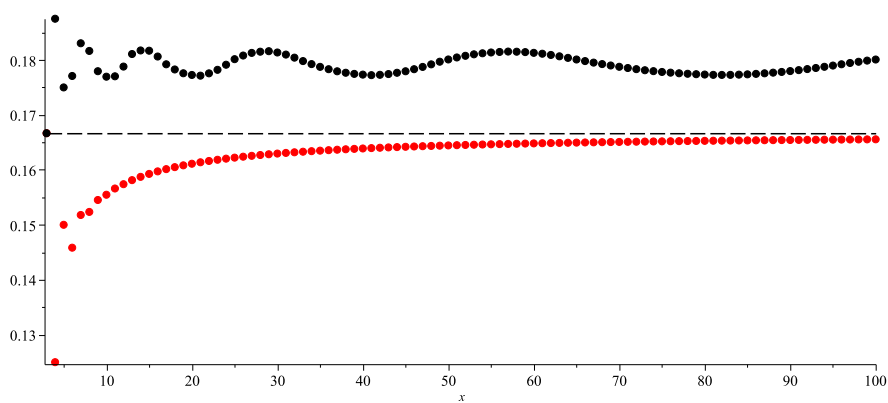}
\caption{The expected proportion of napkinless diners, when seated using trap setting (in red) and napkin shunning (in black). The dashed line is at height $1/6$.}\label{fig: plots}
\end{figure}

Note that the proportion for trap setting has a limit of $1/6$ as $n\to \infty$, agreeing with Winkler's comments and our intuitive understanding of trap setting. We show the plots of these two proportions for tables of size $n=3,4,\ldots,100$, in Figure~\ref{fig: plots}. At a table for 48 diners as in the original puzzle statement, trap setting yields an expected $7.8889$ napkinless diners, while napkin shunning has an expectation of $8.6015$ napkinless diners.

The rest of the paper is organized as follows. In Section~\ref{sec:winkler} we study Winkler's trap setting algorithm and use exact enumerative results to find the generating function (and an exact formula) for the expected number of napkinless diners. In Section~\ref{sec:shun} we study napkin shunning, and although we cannot give an exact formula for the expected number of napkinless diners, we do give bounds that show it is superior to trap setting. In Section \ref{sec:conclusion} we mention several unresolved questions. In particular, we present a variation on napkin shunning which experimentally outperforms the algorithm of Section \ref{sec:shun}, and we make no claims about an optimal strategy for the adaptive maitre d'.

\section{Trap setting}\label{sec:winkler}

The input for any adaptive strategy is a sequence of diner preferences for their right napkin or left napkin, which we call a \textbf{preference order}. We model a preference order for $n$ diners with a list $\sigma = (\sigma_1, \sigma_2, \ldots, \sigma_n)$, where $\sigma_j = -1$ means the $j$th diner prefers the napkin to their left and $\sigma_j = +1$ means the diner prefers the napkin to their right. That is, $\sigma \in \{-1,1\}^n$.

If there are $n$ seats at the table we label positions from left to right from the perspective of the diners (counterclockwise if viewing the table from above) and the diners are labeled by their order of arrival at the table. We choose to declare seat ``1'' as the seat in which the first diner sits. Thus as we move to the right from diner 1, we encounter seats 2, 3, and so on. In moving to the left of diner 1 we find seats $n$, $n-1$, and so on. We write seating orders as lists $w=(w_1,w_2,\ldots,w_n)$, where $w_i = j$ means that the $j$th diner sits in seat $i$. In principle, the maitre d' could create a seating order for any permutation $w$ of $\{1,2,\ldots,n\}$ such that $w_1=1$.

\begin{figure}
\[
\begin{tikzpicture}[scale=1,baseline=0]
\draw (0,0) circle (2);
\draw[thick] (360/8*0+90:2.25) node[fill=white,inner sep=0] {$1$} -- (360/8*0+112.5:1.75) node[circle,fill=white,inner sep=1] {$\curlywedge$};
\draw (360/8*1+90:2.25) node[fill=white,inner sep=0] {$\color{red}{\textbf{5}}$};
\draw[thick] (360/8*2+90:2.25) node[fill=white,inner sep=0] {$2$} -- (360/8*1+112.5:1.75) node[circle,fill=white,inner sep=1] {$\curlywedge$};
\draw[thick] (360/8*3+90:2.25) node[fill=white,inner sep=0] {$8$} -- (360/8*2+112.5:1.75) node[circle,fill=white,inner sep=1] {$\curlywedge$};
\draw[thick] (360/8*4+90:2.25) node[fill=white,inner sep=0] {$4$} -- (360/8*4+112.5:1.75) node[circle,fill=white,inner sep=1] {$\curlywedge$};
\draw[thick] (360/8*5+90:2.25) node[fill=white,inner sep=0] {$6$} -- (360/8*5+112.5:1.75) node[circle,fill=white,inner sep=1] {$\curlywedge$};
\draw (360/8*6+90:2.25) node[fill=white,inner sep=0] {$\color{red}{\textbf{7}}$};
\draw[thick] (360/8*7+90:2.25) node[fill=white,inner sep=0] {$3$} -- (360/8*6+112.5:1.75) node[circle,fill=white,inner sep=1] {$\curlywedge$};
\draw (360/8*3+112.5:1.75) node {$\curlywedge$};
\draw (360/8*7+112.5:1.75) node {$\curlywedge$};
\end{tikzpicture} \quad \longleftrightarrow \quad \begin{array}{cccccccc}
   R&R&L&L&R&L&R&L \\
   \hline
    1&\color{red}{\textbf{5}}&2&8&4&6&\color{red}{\textbf{7}}&3\\
 \end{array}
\]
\caption{A seating arrangement for eight, with preferences that lead to two napkinless diners.}\label{fig:circle}
\end{figure}
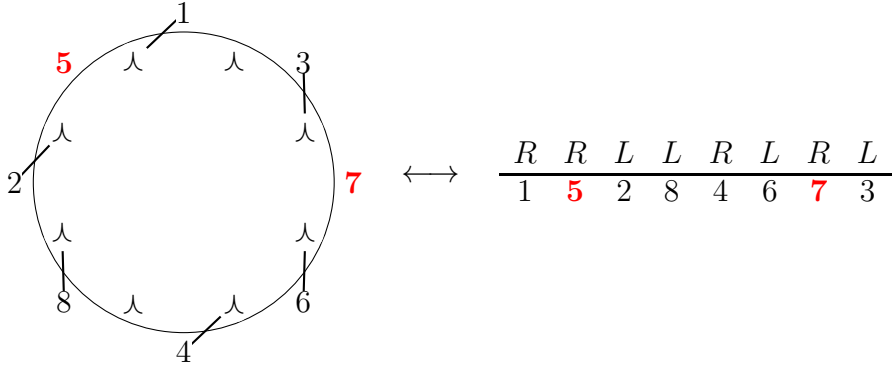

When providing examples, we will often draw tables as arrays with the seating order $w$ on bottom and the corresponding diner preferences above (though usually with $R$ and $L$ instead of $+1$ and $-1$). In Figure \ref{fig:circle} we see two visualizations of the seating order $w=(1,5,2,8,4,6,7,3)$ together with preference order $\sigma=(1,-1,-1,1,1,-1,1,-1)$. Notice that the sequence of events here conspire to leave two of the diners (diners 5 and 7) without a napkin.

\subsection{Algorithm $W$ (the trap setting maitre d')}\label{sec:walgorithm}

We now describe an algorithm for implementing Winkler's trap setting strategy. 

We initialize by seating the first diner and observing which napkin they select. We call this first choice of left or right the \textbf{primary direction}. All movement then proceeds in this primary direction. For simplicity below, we assume the primary is to the right.
\begin{enumerate}
    \item[\textbf{W1.}] (Gap check) If the two seats immediately to the right of the previous diner are empty, go to step W2. Else, go to step W4.

    \item[\textbf{W2.}] (Trap setting) If the previous diner reached right, seat the next diner two seats to the right, and return to step W1. Else, go to step W3.
    
    \item[\textbf{W3.}] (Reset) Seat the next diner one seat to the right of the previous diner.  Return to step W1.
    
    \item[\textbf{W4.}] (Trap springing) Seat remaining diners in the open seats, moving right from the first diner.\end{enumerate}

We will use the notation $\nu_W(\sigma)$ to denote the number of napkinless diners that result from applying algorithm $W$ to preference order $\sigma$.\footnote{In \cite{CP}, the authors approached Conway's napkin problem by viewing the number of napkinless diners as a statistic for signed permutations, i.e., pairs $(w,\sigma)$. In this paper, for a fixed strategy, the number of napkinless diners depends only on $\sigma$.} For example, suppose the preference order for a group of 18 diners is 
\[
 \sigma = (1,1,-1,1,-1,1,1,-1,1,-1,-1,1,1,-1,-1,1,1,-1).
\]

\begin{figure}
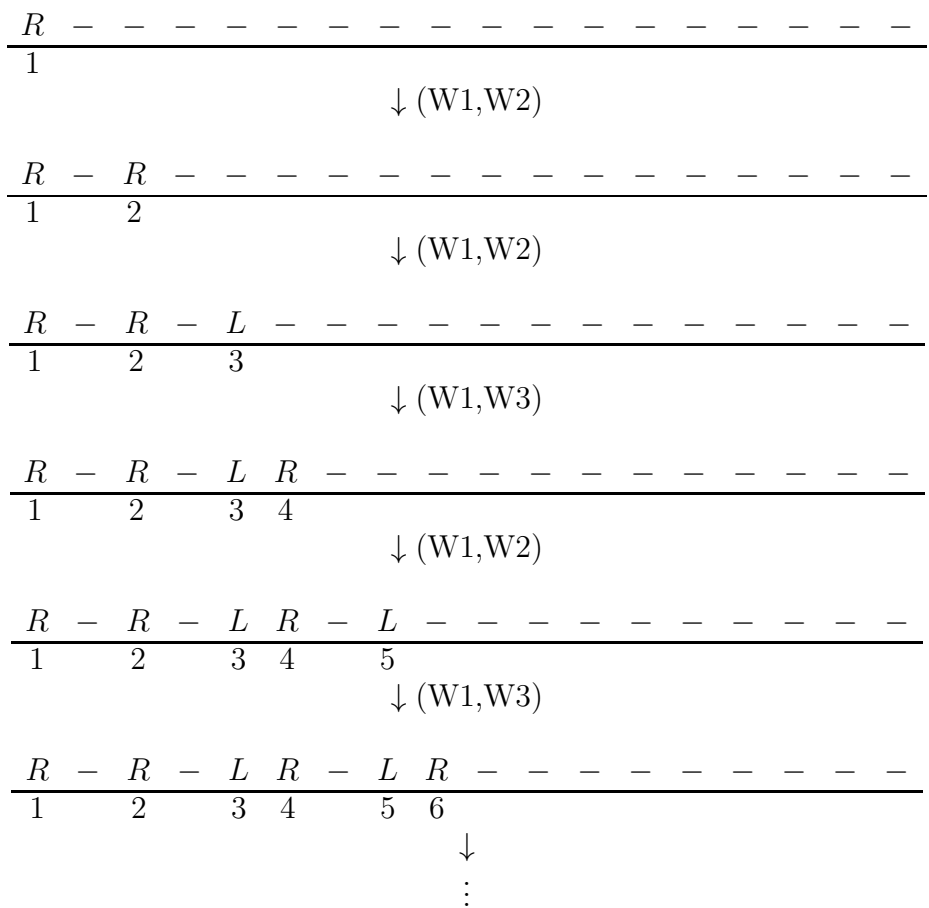

\[
\begin{array}{c}
  \begin{array}{cccccccccccccccccc}
   R&-&-&-&-&-&-&-&-&-&-&-&-&-&-&-&-&- \\
   \hline
   1&&&&&&&&&&&&&&&&&\\
 \end{array}\\
  \downarrow (\text{W1,W2}) \\\\
 \begin{array}{cccccccccccccccccc}
   R&-&R&-&-&-&-&-&-&-&-&-&-&-&-&-&-&- \\
   \hline
   1&&2&&&&&&&&&&&&&&&\\
 \end{array}\\
 \downarrow (\text{W1,W2})\\\\
 \begin{array}{cccccccccccccccccc}
   R&-&R&-&L&-&-&-&-&-&-&-&-&-&-&-&-&- \\
   \hline
   1&&2&&3&&&&&&&&&&&&&\\
 \end{array}\\
 \downarrow (\text{W1,W3})\\\\
 \begin{array}{cccccccccccccccccc}
   R&-&R&-&L&R&-&-&-&-&-&-&-&-&-&-&-&- \\
   \hline
   1&&2&&3&4&&&&&&&&&&&&\\
 \end{array}\\
 \downarrow (\text{W1,W2}) \\\\
 \begin{array}{cccccccccccccccccc}
   R&-&R&-&L&R&-&L&-&-&-&-&-&-&-&-&-&- \\
   \hline
   1&&2&&3&4&&5&&&&&&&&&&\\
 \end{array}\\
 \downarrow (\text{W1,W3})\\\\
 \begin{array}{cccccccccccccccccc}
   R&-&R&-&L&R&-&L&R&-&-&-&-&-&-&-&-&- \\
   \hline
   1&&2&&3&4&&5&6&&&&&&&&&\\
 \end{array}\\
 \downarrow \\
 \vdots\\
\end{array}
\]
\caption{Trap setting for the first six diners for the example preference order.}\label{fig:ex1}
\end{figure}

In Figure \ref{fig:ex1}, we see an illustration of algorithm $W$ applied to this preference order up to the first six diners. At this stage of the seating process, the maitre d' has already set two traps, one between diners $2$ and $3$, and another between diners $4$ and $5$. Continuing through the algorithm, after the first 11 diners have been seated, the table appears as follows:
\[
\begin{array}{cccccccccccccccccc}
   R&-&R&-&L&R&-&L&R&-&R&-&L&R&-&L&L&- \\
   \hline
    1&&2&&3&4&&5&6&&7&&8&9&&10&11\\
 \end{array} .
\]

Notice that after seating these 11 diners, the maitre d' has set four traps, and there is no longer a gap of two open seats after the most recently seated diner.  The algorithm thus proceeds to W4, in which the remaining diners are sequentially seated in the  open seats.  Some of these diners are doomed to spring the four napkinless traps already set and we obtain the following seating arrangement:
\[
\begin{array}{cccccccccccccccccc}
   R&R&R&R&L&R&L&L&R&L&R&R&L&R&R&L&L&L \\
   \hline
    1&12&2&\color{red}{\textbf{13}}&3&4&\color{red}{\textbf{14}}&5&6&15&7&\color{red}{\textbf{16}}&8&9&\color{red}{\textbf{17}}&10&11&18\\
 \end{array}\, .
\]

Thus diners $13$, $14$, $16$, and $17$ are napkinless, and $\nu_W(\sigma) = 4$ for this example. 

\subsection{Table types}

Let $\C_n$ denote a circular table with $n$ empty seats and $n$ napkins between the seats, prior to any diners being seated. We make the same observation as in \cite[Section 2]{CP}, and notice that once the first diner has been seated, circular symmetry is broken, and the dynamics of the rest of the seating takes place on one of four types of linear segments. That is, our seating diagram will appear as concatenated blocks of the following types:
\[
R-\cdots -R, \quad L-\cdots -L, \quad R-\cdots-L, \quad L-\cdots-R,
\]
as in our example of the previous subsection. To aid in our later analysis, we classify these below, where each table type is indexed by the total number of empty seats, $n$.

\begin{itemize}
\item \textbf{The right-leaning table.} This table, denoted $\R_n$, is a table with $n$ empty seats and $n$ napkins, such that the leftmost seat has no napkin to its left.
\item \textbf{The left-leaning table.} This table, denoted $\Lf_n$, is a table with $n$ empty seats and $n$ napkins, such that the rightmost seat has no napkin to its right.
\item \textbf{The inner-facing table.} This table, denoted $\I_n$, is a table with $n$ empty seats and $n-1$ napkins, such that the leftmost seat has no napkin to its left and the rightmost seat has no napkin to its right.
\item \textbf{The outer-facing table.} This table, denoted $\Ou_n$, is a table with $n$ empty seats and $n+1$ napkins, such that each seat has a napkin to both its left and its right.
\end{itemize}

For a given preference order $\sigma$, the number of napkinless diners on a circular table might or might not agree with the number of napkinless diners when those same diners are seated on one of the linear table types. For example, if $\sigma_i =+1$ for all $i$, everyone will get a napkin on table $\C_n$, $\R_n$, or $\Ou_n$, but the pigeonhole principle dictates that at least one person will always be napkinless on table $\I_n$. To help distinguish, we write $\nu_W(\mathcal{T}_n;\sigma)$ to make clear that algorithm $W$ is being used to seat $n$ diners with preference order $\sigma$ on table type $\mathcal{T}_n$. Further, we will use the notation $W(\mathcal{T}_n;t)$ to denote the generating function for the distribution of napkinless numbers over all preference orders, i.e.,
\[
 W(\mathcal{T}_n;t) = \sum_{\sigma \in \{-1,1\}^n} t^{\nu_W(\mathcal{T}_n;\sigma)}.
\]

\subsection{Generating functions}\label{sec:Wgf}

For each circular table of size $n$, let $W_n(t)=W(\C_n;t)$ be the generating function for the number of napkinless at a table of size $n$ over all preference orders. That is,
\[
 W_n(t) = \sum_{\sigma \in \{-1,1\}^n} t^{\nu_W(\C_n ;\sigma)} = \sum_{k=0}^{\lfloor n/3\rfloor} a_{n,k}\cdot t^k.
\]
Here $a_{n,k}$ is the number of preference orders that result in $k$ napkinless diners when seated according to algorithm $W$:
\[
a_{n,k} = \left|\left\{ \sigma \in \{-1,1\}^n : \nu_W(\C_n;\sigma) = k\right\}\right|.
\]
For example, $W_3(t) = 4+4t$, $W_4(t) = 8+8t$, and $W_5(t) = 8+24t$.

We can notice that for $n\geq 2$, the numbers $a_{n,k}$ are all multiples of 4. This follows because $\nu_W(\C_n;-\sigma)=\nu_W(\C_n;\sigma)=\nu_W(\C_n;\sigma')=\nu_W(\C_n;-\sigma')$, where $\sigma'$ is the preference order whose entries are identical to $\sigma$ except in the final position. That is, both the action of toggling the preference of the final diner and the action of reversing all preferences leave the number of napkinless diners unchanged. We let $b_{n,k} = a_{n,k}/4$, and we display the numbers $b_{n,k}$ for values of $n\leq 14$ in Table \ref{tab:bnk}. Also provided in the table is the expected number of napkinless diners when seated with algorithm $W$.

\begin{table}
\[
\begin{array}{c|ccccc||c}
n\backslash k & 0 & 1 & 2 & 3 & 4 & E_n[\nu_W(\sigma)] \\
\hline 
\hline
2 & 1 &&&&& 0\\
3 & 1 & 1 &&&& 1/2\\
4 & 2 & 2 &&&& 1/2\\
5 & 2 & 6 &&&& 3/4\\
6 & 4 & 10 & 2 &&& 7/8\\
7 & 4 & 22 & 6 &&& 17/16\\
8 & 8 & 34 & 22 &&& 39/32\\
9 & 8 & 66 & 50 & 4 && 89/64\\
10 & 16 & 98 & 126 & 16 && 199/128\\
11 & 16 & 178 & 250 & 68 && 441/256 \\
12 & 32 & 258 & 534 & 192 & 8 & 967/512\\
13 & 32 & 450 & 978 & 548 & 40 & 2105/1024\\
14 & 64 & 642 & 1902 & 1296 & 192 & 4551/2048
\end{array}
\]
\caption{The numbers $b_{n,k}$, such that $4b_{n,k}$ is the number of preference orders of length $n$ that result in $k$ napkinless diners using algorithm $W$. In the final column, we show the expected number of napkinless diners for a circular table of size $n$.}\label{tab:bnk}
\end{table}

\subsection{Recursive structure}\label{sec:recur}

Upon seating the first diner, the remainder of the algorithm plays out on either a table of type $\R_{n-1}$ or of type $\Lf_{n-1}$, depending on whether the primary direction $\sigma_1$ is $+1$ or $-1$. Thus, 
\[
 W_n(t) = W(\C_n;t) = W(\R_{n-1};t) + W(\Lf_{n-1};t).
\]
By the mirror symmetry noted previously, $\nu_W(\R_n;\sigma) = \nu_W(\Lf_n;-\sigma)$, so the two cases have the same distributions of napkinless diners: $W(\R_{n-1};t)=W(\Lf_{n-1};t)$. Hence
\begin{equation}\label{eq:Cnid}
 W_n(t) = 2W(\R_{n-1};t).
\end{equation}

Now we want to consider what happens with a right-leaning table with primary direction $R$. Then the next step of the  algorithm is as follows. In the images below, diner 1 is indicated twice to highlight the presence of the linear table type. 
\[
\begin{array}{cccccccccccccc}
   R&-&-&-&\cdots &-&-&R\\
   \hline
   1&&&&&&&1\\
 \end{array}
  \longrightarrow 
 \begin{array}{c}
 \begin{array}{cccccccccccccc}
   R&-&R&-&\cdots &-&-&R \\
   \hline
   1&&2&&&&&1\\
 \end{array} \\
 \text{or} \\
 
  \begin{array}{cccccccccccccc}
   R&-&L&-&\cdots &-&-&R \\
   \hline
   1&&2&&&&&1\\
 \end{array} .
 \end{array}
\]
In terms of our generating function notation, this shows
\[
W(\R_n;t) = W(\R_1;t)W(\R_{n-2};t)+W(\I_1;t)W(\Ou_{n-2};t),
\]
where $n$ here is the number of empty seats in the block on the left. 

As $W(\R_1;t) = 2$ and $W(\I_1;t) = 2t$, we have
\begin{equation}\label{eq:Rn}
W(\R_n;t) = 2W(\R_{n-2};t)+2tW(\Ou_{n-2};t).
\end{equation}

Now, let's observe how the algorithm treats a block of type $\Ou_n$. As before, without loss of generality, we can suppose the primary direction is right, and hence the next diner sits in the leftmost empty seat of the table segment below. Algorithm $W$ then proceeds as:
\[
\begin{array}{cccccccccccccc}
   L&-&-&-&\cdots &-&-&R\\
   \hline
  2&&&&&&&1\\
 \end{array}
  \longrightarrow 
 \begin{array}{c}
 \begin{array}{cccccccccccccc}
   L&R&-&-&\cdots &-&-&R\\
   \hline
  2&3&&&&&&1\\
 \end{array} \\
 \mbox{ or }  \\
 \begin{array}{cccccccccccccc}
   L&L&-&-&\cdots &-&-&R\\
   \hline
  2&3&&&&&&1\\
 \end{array}.
 \end{array}
\]
In terms of our generating functions,  we see that 
\begin{equation}\label{eq:Orec}
 W(\Ou_n;t) = W(\R_{n-1};t)+W(\Ou_{n-1};t).
\end{equation}

Now we can suppose that $n\geq 3$ and judiciously apply \eqref{eq:Rn} and \eqref{eq:Orec} to find:
\begin{align*}
W(\R_n;t) &=2W(\R_{n-2};t)+2tW(\Ou_{n-2};t), \\
 &=2W(\R_{n-2};t) + 2t( W(\R_{n-3};t)+W(\Ou_{n-3};t)), \\
 &=2W(\R_{n-2};t) + 2tW(\R_{n-3};t)+2tW(\Ou_{n-3};t),\\
 &=2W(\R_{n-2};t) + 2tW(\R_{n-3};t)+W(\R_{n-1};t)-2W(\R_{n-3};t),\\
 &=W(\R_{n-1};t)+2W(\R_{n-2};t)+2(t-1)W(\R_{n-3};t).
\end{align*} 
Multiplying by 2 in consideration of Equation \eqref{eq:Cnid} yields the recurrence in the following proposition.

\begin{proposition}\label{prp:Wrec}
For $n\geq 4$, we have 
\begin{equation}\label{eq:Wrec}
 W_n(t) = W_{n-1}(t) + 2W_{n-2}(t) + 2(t-1)W_{n-3}(t).
\end{equation}
\end{proposition}

By setting $W_0(t)=2$, we can ensure the above identity holds for all $n\geq 3$. Now define the ordinary generating function for the polynomials $W_n(t)$ as
\begin{align*}
W(t,z) &= \sum_{n\geq 0} W_n(t)z^n \\
 &= 2+2z+4z^2+(4+4t)z^3 + (8+8t)z^4 + (8+24t)z^5 + \cdots.
\end{align*}
With a dash of ``generatingfunctionology'' (see \cite{Wilf}) we can use the recurrence in Equation \eqref{eq:Wrec} to manipulate this series and obtain the following rational expression for $W(t,z)$ as a corollary to Proposition \ref{prp:Wrec}.

\begin{corollary}\label{cor:Wgf}
We have
\begin{equation}\label{eq:Wgf}
W(t,z) = \frac{2-2z^2}{1-z-2z^2-2(t-1)z^3}.
\end{equation}
\end{corollary}

\subsection{Expectations for trap setting}

The expected value of $\nu_W(\C_n;\sigma)$, which we denote $E_n^W=E_n[\nu_W(\C_n;\sigma)]$, is easily computed from $W_n(t)$ via
\[
 E_n^W = \sum_{k=0}^{\lfloor n/3 \rfloor} k\cdot \frac{a_{n,k}}{2^n} = \frac{W_n'(1)}{2^n}.
\]
We provide these values in the rightmost column of Table \ref{tab:bnk} for $n\leq 14$.

By taking the derivative at $t=1$ on both sides of \eqref{eq:Wrec} we obtain the following recurrence for all $n\geq 4$:
\begin{equation}\label{eq:Erec}
 E_n^W = \frac{E^W_{n-1}}{2} + \frac{E^W_{n-2}}{2} + \frac{1}{4}.
\end{equation}

Now it is another exercise in generating functions to obtain the following result, which proves part (1) of Theorem \ref{thm:main}.

\begin{proposition}\label{prp:EW}
The generating function for the expected number of napkinless diners with algorithm $W$ is:
\[
 E(z) = \sum_{n\geq 0} E^W_n z^n = \frac{z^3(2-z)}{2(1-z)^2(2+z)}.
\]
For $n\geq 0$, the formula for the expectation on a table for $n+3$ diners is
\[
E^W_{n+3}=\frac{(3n+7)2^{n-1}+(-1)^n}{9\cdot 2^n}.
\]
\end{proposition}

The rational expression for generating function $E(z)$ follows either by building up directly from the recurrence in \eqref{eq:Erec}, or by differentiating \eqref{eq:Wgf} with respect to $t$ and substituting $z/2$ for $z$. Obtaining the formula for the expectation is a straightforward process for rational generating functions such as this one, with only linear factors in the denominator. 

We remark that the sequence $2^nE^W_{n+3}$, which, for $n\geq 1$ begins
\[
1, 3, 7, 17, 39, 89, 199, 441, 967, \ldots
\] 
is precisely entry A127984 of the OEIS \cite{oeis}.

\section{Napkin shunning}\label{sec:shun}

In this section we study another adaptive algorithm for the malicious maitre d' that performs strictly better than algorithm $W$. As with Algorithm $W$, our input is a preference order $\sigma$, and the maitre d' makes a choice about where to seat diner $i$ based on the observed preferences of diners $1, \ldots, i-1$.

\subsection{Algorithm $S$ (the napkin shunning maitre d')} We initialize by seating the first diner and observing which napkin they select. The algorithm then proceeds as follows.
\begin{enumerate}
\item[\textbf{S1.}] (Gap check) If all seats are filled, stop. If the seat adjacent to the previous diner in the direction opposite their preference is open, go to step S2. Else, go to step S3.
\item[\textbf{S2.}] (Napkin shunning) Place the next diner in the seat identified in step S1. Return to step S1. 
\item[\textbf{S3.}] (Reset) Place the new diner $\lceil i/2 \rceil$ seats from the left edge of the leftmost gap, where $i$ is the number of empty seats in this gap. Return to step S1. 
\end{enumerate}

We write $\nu_S(\sigma)$ to denote the number of napkinless diners that result from applying algorithm $S$ to preference order $\sigma$. For example, suppose the preferences of 18 diners are given by the same preference sequence as we used in Section \ref{sec:walgorithm}:
\[
 \sigma = (1,1,-1,1,-1,1,1,-1,1,-1,-1,1,1,-1,-1,1,1,-1).
\]

\begin{figure}
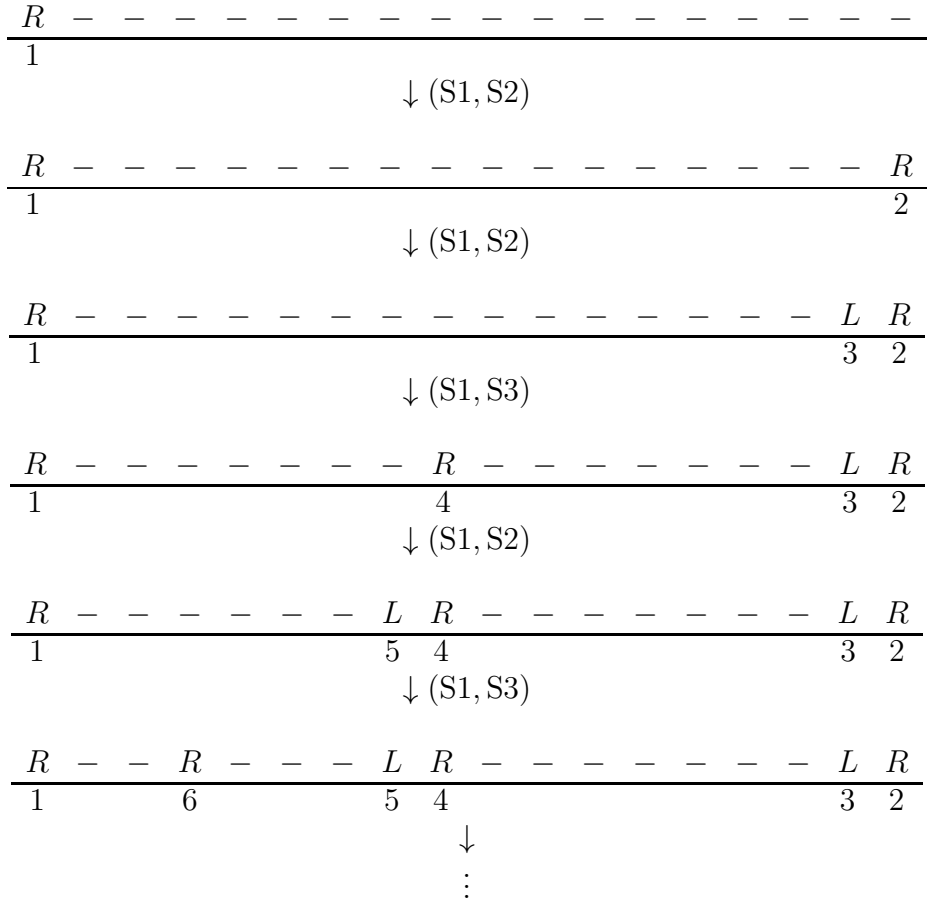

\[
\begin{array}{c}
  \begin{array}{cccccccccccccccccc}
   R&-&-&-&-&-&-&-&-&-&-&-&-&-&-&-&-&- \\
   \hline
   1&&&&&&&&&&&&&&&&&\\
 \end{array}\\
  \downarrow  (\text{S1},\text{S2}) \\\\
 \begin{array}{cccccccccccccccccc}
   R&-&-&-&-&-&-&-&-&-&-&-&-&-&-&-&-&R \\
   \hline
   1&&&&&&&&&&&&&&&&&2\\
 \end{array}\\
 \downarrow (\text{S1},\text{S2}) \\\\
 \begin{array}{cccccccccccccccccc}
   R&-&-&-&-&-&-&-&-&-&-&-&-&-&-&-&L&R \\
   \hline
   1&&&&&&&&&&&&&&&&3&2\\
 \end{array}\\
 \downarrow (\text{S1},\text{S3}) \\\\
 \begin{array}{cccccccccccccccccc}
   R&-&-&-&-&-&-&-&R&-&-&-&-&-&-&-&L&R \\
   \hline
   1&&&&&&&&4&&&&&&&&3&2\\
 \end{array}\\
 \downarrow (\text{S1},\text{S2}) \\\\
\begin{array}{cccccccccccccccccc}
   R&-&-&-&-&-&-&L&R&-&-&-&-&-&-&-&L&R \\
   \hline
   1&&&&&&&5&4&&&&&&&&3&2\\
 \end{array}\\
 \downarrow (\text{S1},\text{S3})\\\\
 \begin{array}{cccccccccccccccccc}
   R&-&-&R&-&-&-&L&R&-&-&-&-&-&-&-&L&R \\
   \hline
   1&&&6&&&&5&4&&&&&&&&3&2\\
 \end{array}\\
 \downarrow \\
 \vdots
\end{array}
\]
\caption{Napkin shunning for the first six diners for the example preference order.}\label{fig:ex2}
\end{figure}

In Figure \ref{fig:ex2}, we see an illustration of algorithm $S$ applied to this preference order up to the first six diners. At this point, we have already shunned two napkins, one between diners $2$ and $3$, and another between diners $4$ and $5$. Diners $7$ through $11$ will proceed to fill rest of the seats between diners $1$ and $4$, until the table appears as follows:
\[
\begin{array}{cccccccccccccccccc}
   R&L&R&R&L&R&L&L&R&-&-&-&-&-&-&-&L&R \\
   \hline
   1&8&7&6&10&9&11&5&4&&&&&&&&3&2\\
 \end{array} \, .
\]
Notice that diner $4$ was the first to be seated as a result of step S3. The algorithm now proceeds to fill in the gap to the right of diner $4$. The final seating arrangement is:
\[
\begin{array}{cccccccccccccccccc}
   R&L&R&R&L&R&L&L&R&L&L&R&R&R&R&L&L&R \\
   \hline
   1&8&7&6&10&9&\color{red}{\textbf{11}}&5&4&\color{red}{\textbf{15}}&14&13&12&17&16&\color{red}{\textbf{18}}&3&2\\
 \end{array}\, ,
\]
and we can see that diners $11$, $15$, and $18$ get stuck without a napkin. 

It is worth noting that for this particular preference order, Algorithm $W$ produced four napkinless diners but Algorithm $S$ only produced three, i.e., for this $\sigma$, $\nu_W(\sigma) = 4 > 3 = \nu_S(\sigma)$. In other words, the napkin shunning algorithm is not universally better than trap setting. Nonetheless, we will see that its expected number of napkinless diners taken over all sign sequences is greater.

\subsection{Generating functions}

We define the generating function for algorithm $S$ just as we did with algorithm $W$, tweaking the notation only slightly. Let
\[
 S_n(t) = \sum_{\sigma \in \{-1,1\}^n} t^{\nu_S(\C_n;\sigma)}= \sum_{k=0}^{\lfloor n/3\rfloor} c_{n,k} t^k,
\]
so that $c_{n,k}$ is the number of sign sequences that result in $k$ napkinless diners when seated according to algorithm $S$:
\[
c_{n,k} = \left|\left\{ \sigma \in \{-1,1\}^n : \nu_S(\C_n;\sigma) = k\right\}\right|.
\]
For example, $S_3(t) = 4+4t$ and $S_4(t) = 4+12t$, and $S_5(t) = 4+28t$. We can see the numbers $c_{n,k}$ are multiples of 4 for exactly the same reasons the numbers $a_{n,k}$ were multiples of 4, as discussed in Section \ref{sec:Wgf}. In a similar fashion, we define the numbers $d_{n,k}=c_{n,k}/4$. We see the $d_{n,k}$ for values of $n\leq 14$ in Table \ref{tab:dnk}, along with the expected number of napkinless diners when seated with algorithm $S$. 

\begin{table}
\[
\begin{array}{c|ccccc||c}
n\backslash k & 0 & 1 & 2 & 3 & 4 & E_n[\nu_S(\sigma)] \\
\hline 
\hline
2 & 1 &&&&& 0\\
3 & 1 & 1 &&&& 1/2\\
4 & 1 & 3 &&&& 3/4\\
5 & 1 & 7 &&&& 7/8\\
6 & 1 & 13 & 2 &&& 17/16\\
7 & 1 & 21 & 10 &&& 41/32\\
8 & 1 & 33 & 30 &&& 93/64\\
9 & 1 & 49 & 78 &  && 205/128\\
10 & 1 & 69 & 174 & 12 && 453/256\\
11 & 1 & 93 & 350 & 68 && 997/512 \\
12 & 1 & 121 & 638 & 256 & 8 & 2197/1024\\
13 & 1 & 153 & 1086 & 736 & 72 & 4821/2048\\
14 & 1 & 193 & 1790 & 1800 & 312 & 10421/4096
\end{array}
\]
\caption{The numbers $d_{n,k}$, such that $4d_{n,k}$ is the number of sign sequences of length $n$ that result in $k$ napkinless diners using Algorithm $S$. The final column shows the expected number of napkinless diners for a circular table of size $n$.}\label{tab:dnk}
\end{table}

\subsection{Recursive structure, revisited}

We will revisit the pictographic understanding of the seating algorithm as we did in Section \ref{sec:recur}, but now for the napkin shunning algorithm.

The initial step is the same. After the first diner has been seated, the rest of the algorithm plays out on either $\R_{n-1}$  or on $\Lf_{n-1}$. Thus,
\[
 S_n(t) = S(\C_n;t) = S(\R_{n-1};t) + S(\Lf_{n-1};t),
\]
and by the mirror symmetry $\nu_S(\R_n;\sigma) = \nu_S(\Lf_n;-\sigma)$, we again have
\begin{equation}\label{eq:Cnid2}
 S_n(t) = 2S(\R_{n-1};t).
\end{equation}

Now we want to consider what happens with a right-leaning table with key direction $R$. Then the next step of algorithm $S$ looks like:
\[
\begin{array}{cccccccccccccc}
   R&-&-&\cdots &-&-&-&R\\
   \hline
   1&&&&&&&1\\
 \end{array}
 \longrightarrow 
 \begin{array}{c}
 \begin{array}{cccccccccccccc}
   R&-&-&\cdots &-&-&R&R \\
   \hline
   1&&&&&&2&1\\
 \end{array} \\
 \text{or} \\
  \begin{array}{cccccccccccccc}
   R&-&-&\cdots &-&-&L&R \\
   \hline
   1&&&&&&2&1\\
 \end{array},
 \end{array}
\]
which in terms of generating functions shows
\begin{equation}\label{eq:Rn2}
S(\R_n;t) = S(\R_{n-1};t)+S(\I_{n-1};t).
\end{equation}

Now we study algorithm $S$ on the inner-leaning table $\I_n$. Here 
 algorithm $S$ proceeds as:
\[
\begin{array}{cccccccccccccc}
   R&-&-&\cdots &-&-&L \\
   \hline
   1&&&&&&2\\
 \end{array}
 \longrightarrow
 \begin{array}{c}
 \begin{array}{cccccccccccccc}
   R&-&-&\cdots & R & \cdots &-&-&L \\
   \hline
   1&&&&3&&&&2\\
 \end{array} \\
 \mbox{ or }\\
  \begin{array}{cccccccccccccc}
   R&-&-&\cdots & L & \cdots &-&-&L \\
   \hline
   1&&&&3&&&&2\\
 \end{array},
 \end{array}
\]
where, if there are $n$ empty seats, diner 3 has been placed in seat $\lceil n/2 \rceil$. The number of empty seats to the left of diner 3 is $\lfloor (n-1)/2 \rfloor$ and the number of empty seats to the right is $\lceil (n-1)/2 \rceil$. Thus, 
\begin{equation}\label{eq:Irec}
 S(\I_n;t) = S(\R_{\lfloor \frac{n-1}{2}\rfloor}; t) S(\I_{\lceil \frac{n-1}{2} \rceil}; t) + S(\I_{\lfloor \frac{n-1}{2}\rfloor}; t) S(\R_{\lceil \frac{n-1}{2}\rceil}; t),
\end{equation}
where we have used the identity $S(\R_k;t) = S(\Lf_k;t)$ in the second term.

Let us pause for a moment now and consider the nature of the recurrences we have developed in Equations \eqref{eq:Cnid2}, \eqref{eq:Rn2}, and \eqref{eq:Irec}. These three intertwined identities can be used, with the help of a computer, to recursively compute $S_n(\C_n;t) = S_n(t)$ for the first few hundred values of $n$. However, unlike the case for the $W_n(t)$, we do not have a simple way to ``solve for'' $S_n(t)$ in terms of itself with a smaller index. Using \eqref{eq:Rn2} repeatedly, we can establish the relationship
\begin{align*}
 S(\R_n;t) &= S(\R_{n-1};t) + S(\I_{n-1};t),\\
  &=S(\R_{n-2};t) + S(\I_{n-2};t) + S(\I_{n-1};t),\\
  &\vdots \\
  &=S(\R_1;t) + S(\I_1;t) + \cdots + S(\I_{n-1};t),\\
  &=2 + \sum_{i=1}^{n-1}S(\I_i;t).
\end{align*}
From this identity, we can obtain a functional identity linking the generating function for the polynomials $S_n(t)$ to the generating function for the polynomials $S(\I_n;t)$. However, this begs the question: what is the generating function for the $S(\I_n;t)$?

The only other identity at our disposal is Equation \eqref{eq:Irec}. While at first glance it seems innocent (a two-term quadratic recurrence isn't \emph{so} awful\ldots) we notice that the indices on the right hand side are not a constant distance from the index on the left. That is, this is not a fixed length recurrence. If it were, the generating function would be algebraic, but as it stands, all bets are off. Let us avoid the exact distribution question, then, and proceed to study the expected values $E[\nu_S(\sigma)]$. 

\subsection{Expectations for napkin shunning}

Let $E_n^S =E[\nu_S(\C_n;\sigma)]$ denote the expected number of napkinless diners using Algorithm $S$. As before, we obtain $E_n^S$ by differentiation at $t=1$:
\[
 E_n^S = \sum_{k=0}^{\lfloor n/3 \rfloor} k\cdot \frac{c_{n,k}}{2^n} = \frac{S_n'(1)}{2^n}.
\]
These values are shown in the rightmost column of Table \ref{tab:dnk}. In a similar fashion, define the napkinless expectations for tables $\R_n$ and $\I_n$:
\[
 R_n = E[\nu_S(\R_n;\sigma)] = \frac{\frac{d}{dt}[S(\R_{n};t)]_{t=1}}{2^n}
 \,\mbox{ and }\,  I_n = E[\nu_S^{\I_n}(\sigma)] = \frac{\frac{d}{dt}[S(\I_{n};t)]_{t=1}}{2^n}.
\]

By taking the derivative at $t=1$ on both sides of \eqref{eq:Cnid2}, we get
\begin{equation}\label{eq:ER}
E_n^S = \frac{S_n'(1)}{2^n} = \frac{2\frac{d}{dt}[S(\R_{n-1};t)]_{t=1}}{2^n} = \frac{\frac{d}{dt}[S(\R_{n-1};t)]_{t=1}}{2^{n-1}} = R_{n-1}. 
\end{equation}
Similarly, from 
\eqref{eq:Rn2}, we get
\begin{equation}\label{eq:RI}
 R_n = \frac{1}{2}( R_{n-1} + I_{n-1}),
\end{equation}
and \eqref{eq:Irec} gives
\begin{align*}
I_n &= \frac{1}{2}\left(R_{\lfloor \frac{n-1}{2} \rfloor} + I_{\lceil \frac{n-1}{2} \rceil} + I_{\lfloor \frac{n-1}{2} \rfloor} +R_{\lceil \frac{n-1}{2} \rceil}\right),\\
&= \frac{1}{2}\left(R_{\lfloor \frac{n-1}{2} \rfloor}+ I_{\lfloor \frac{n-1}{2} \rfloor}\right) + \frac{1}{2}\left(R_{\lceil \frac{n-1}{2} \rceil}+I_{\lceil \frac{n-1}{2} \rceil}\right),\\
&=R_{\lfloor \frac{n+1}{2} \rfloor}+R_{\lceil \frac{n+1}{2} \rceil}.
\end{align*}
Putting this together with \eqref{eq:ER} and \eqref{eq:RI} gives us the following result.

\begin{proposition}\label{prp:Erec2}
The expected number of napkinless diners with algorithm $S$ is given by initial values $E^S_1=E^S_2 = 0$, $E^S_3=1/2$, $E^S_4=3/4$, and for $n\geq 5$, by the recurrence
\[
 E^S_n = \frac{1}{2}\left( E^S_{n-1} + E^S_{\lfloor \frac{n+1}{2} \rfloor} + E^S_{\lceil \frac{n+1}{2} \rceil}\right).
\]
\end{proposition}

Now we would like to understand the behavior of $E^S_n/n$ as $n\to \infty$. First, we have a lemma about bounds.

\begin{lemma}\label{lem:bound}
Fix $k\geq 1$. Suppose there exist constants $\alpha$ and $\beta$ such that 
\[
 \alpha \leq \frac{E^S_i}{i} \leq \beta
\]
for each $i=k,\ldots, 2k-2$. Then
\[
 \alpha \leq \frac{E^S_n}{n} \leq \beta
\]
for all $n\geq k$.
\end{lemma}

\begin{proof}
The argument is symmetrical, so we will prove only the lower bound. The induction step is to show that if the bound holds for $i=k,\ldots, 2k-2$, then it holds for $i=k+1,\ldots,2k$. This follows directly from the recurrence in Proposition \ref{prp:Erec2} as we now demonstrate. 

First, consider $i=2k-1$. We have
\begin{align*}
 E^S_{2k-1} &= \frac{1}{2}(E_{2k-2}+E_k+E_k)\\
  &\geq \frac{1}{2}( (2k-2)\alpha + k\alpha + k\alpha),\\
  &=\frac{1}{2}(4k-2)\alpha=(2k-1)\alpha,
\end{align*}
as desired. Similarly,
\begin{align*}
 E^S_{2k} &= \frac{1}{2}(E_{2k-1}+E_k+E_{k+1})\\
  &\geq \frac{1}{2}( (2k-1)\alpha + k\alpha + (k+1)\alpha),\\
  &=\frac{1}{2}4k\alpha=2k\alpha.
\end{align*}
The lemma now follows.
\end{proof}

With $E^S_3/3 =1/6$ and $E^S_4/4 = 3/16$, the lemma immediately tells us that 
\[
8/48\leq E^S_n/n \leq 9/48
\] 
for all $n\geq 3$. Because $41/224 = E^S_5/5 > E^S_6/6 > E^S_7/7 > E^S_8/8=453/2560$, we have
\[
0.1769531250 \approx \frac{453}{2560} =\frac{3171}{17920} \leq \frac{E^S_n}{n} \leq \frac{3280}{17920} =\frac{41}{224} \approx 0.1830357143,
\]
for all $n\geq 5$. This gap size is $109/17920 \approx 0.00608$, so we can say the proportion of napkinless diners is about $18\%$, plus or minus about three tenths of a percent. We have now established part (2) of Theorem \ref{thm:main}, and shown definitively that $E^W_n < n/6 \leq E^S_n$ for all $n\geq 3$.

\section{Further thoughts}\label{sec:conclusion}

We finish this article with some remarks and questions for further study.

\begin{itemize}
\item \textbf{Does $\lim_{n\to \infty} E_n^S/n$ exist?}
Let 
\[
\alpha_k = \min\{ E^S_n/n : n=k,\ldots,2k-2\}
\]
 and 
\[
 \beta_k = \max\{E^S_n/n: n=k,\ldots,2k-2\}.
\]
From Lemma \ref{lem:bound}, it follows that $\alpha_k \leq \alpha_{k+1} \leq \beta_{k+1} \leq \beta_{k}$. Since each sequence is monotone and bounded, the limits
\[
 L= \lim_{k\to \infty} \alpha_k \quad \mbox{ and } \quad U = \lim_{k\to \infty} \beta_k,
\]
must exist.
Empirically, by computing $\alpha_k$ and $\beta_k$ up to $k=5000$, we observe that 
\[
 0.1772860948 \leq L \leq U \leq 0.1814947641.
\]
It would be lovely to show $L\neq U$, but we cannot rule out the possibility that $L=U$. In particular, we do not know if the sequence $E_n^S/n$ converges.

\item \textbf{What are the generating functions for napkin shunning?} Because of the relatively simple recursive structure of Algorithm $W$, we were able to find a rational expression for both the generating function $W(t,z) = \sum W_n(t) z^n$ and the generating function for the corresponding expectations. The recurrences we found for Algorithm $S$ were not as simple, but that does not preclude a different approach to the problem that might be used to deduce some closed form expression for either the generating function
\[
 S(t,z) = \sum_{n\geq 0} S_n(t) z^n \quad \mbox{ or } \quad E^S(z) = \sum_{\neq 0} E_n^S z^n.
\] 

\item \textbf{Is there an optimal strategy for the adaptive maitre d'? (If so, what is it?)} While Algorithm $S$ is superior to Algorithm $W$, it does not provide the optimal adaptive strategy. Let Algorithm $\widetilde{S}$ be identical to Algorithm $S$, except instead of step $S3$, we have step $\widetilde{S3}$, which modifies the place in which we ``split the gap'' in one special case. Specifically, for a gap of size $i>1$, if $\lceil i/2 \rceil \equiv 1 \pmod 3$, we choose to place our next diner $\lceil i/2 \rceil-1$ seats from the left rather than $\lceil i/2 \rceil$ seats from the left. This small tweak to the algorithm produces notable differences in the proportions of napkinless diners. See Figure \ref{fig:pl}, which compares $E_n^S/n \leq E_n^{\widetilde{S}}/n$ for $n=3,4,\ldots,100$. Empirically, $E_n^{\widetilde{S}}/n$ is greater than $18\%$ for $n\geq 12$.

\begin{figure}[h]
\[
\includegraphics[width=12cm]{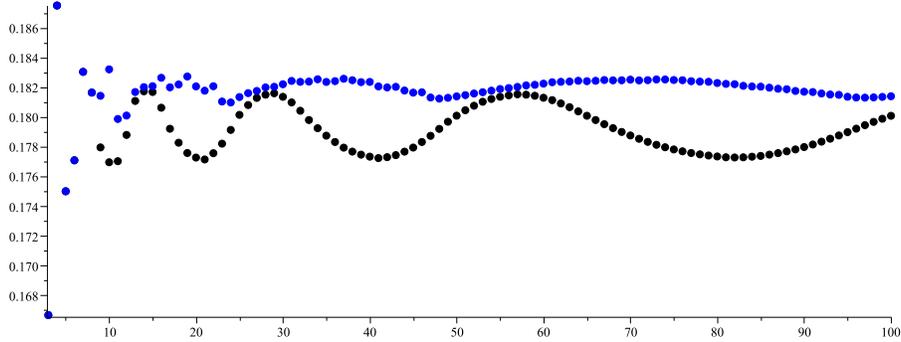} 
\]
\caption{The expected proportion of napkinless diners with the napkin shunning algorithm $S$ (in black) and modified algorithm $\widetilde{S}$ (in blue).}
\label{fig:pl}
\end{figure}

\end{itemize}

There are of course many other variations to the problem of the malicious maitre d' that one could consider. For example, we have implicitly assumed that a diner's preference for left or right is $1/2$, but in \cite{CP, Eriksen, Sudbury} the authors consider the situation of a fixed probability $p$ of taking the left napkin instead. The papers \cite{CP, Sudbury} also consider an extra parameter $f$ that represents the proportion of French diners (who will only take a napkin if it is their preferred napkin) and keep track of more statistics like ``happy diners'' who get the napkin they want and ``frustrated diners'' who get a napkin, but not their preferred napkin. It would be interesting to re-examine the malicious maitre d' puzzle in all its forms with these extra considerations.

\end{document}